\documentclass[11pt,a4paper,english]{amsart}
\usepackage{amssymb,amsmath,amsthm,babel}

\newcommand{\Z}{{\mathbb Z}}
\newcommand{\e}{{\varepsilon}}

\title{A remark on a recent proof of Lehmer's Conjecture.}
\author{Francesco Amoroso}

\date{September 22, 2018}
\begin{document}
\maketitle

\centerline{\it Laboratoire de math\'ematiques Nicolas Oresme, CNRS 
UMR 6139}
\centerline{\it Universit\'e de Caen Normandie, Campus II, BP 5186}
\centerline{\it 14032 Caen Cedex, France}

\section{Introduction}

Let $\beta$ be a non zero algebraic integer and $P_\beta\in\Z[x]$ its minimal polynomial. 
Recall that the Mahler measure $M(\beta)=M(P_\beta)$ is the product of the absolute value of the conjugates of $\beta$ lying outside the unite circle. A classical result of Kronecker ensures that $M(\beta)=1$ if and only if $\beta$ is a root of unity. 
A well-known problem of Lehmer~\cite{Le} asks to find for any positive $\e$ an algebraic integer $\beta$ with $1<M(\beta)<1+\e$. Since then, there is a consensus in a negative answer to Lehmer's question: the so called Lehmer's conjecture becomes the 
existence of an absolute constant $C>1$ such that $M(\beta) \geq C$ for any non zero algebraic integer $\beta$ which is not a root of unity,
Lehmer's Conjecture is still unsolved, the best unconditional result in this direction is a celebrated theorem of Dobrowolski~\cite{Do} which implies that for any $\e>0$ there is a $C(\e)>0$ such that $M(\beta)\geq 1+C(\e)d^{-\e}$ for any non zero algebraic integer $\beta$ of degree $d$ which is not a root of unity.\\

Over one year ago, a very long preprint \cite{JLVG} posted on ArXiv and HAL announced a proof of this conjecture (and of other related results). Unfortunately, as was remarked by several specialists, this proof contains a (at least one) fatal error. The aim of this very short note is to inform the mathematical community, which could be aware of this.

\section{On the proof of Theorem 5.23 of \cite{JLVG}}
The author introduces (definition 4.10 at p. 70) a power series $f_\beta$ which converges on the open unit disc $D(0,1)$. At page 108 he writes ``The key result which makes the link between $M(\beta)$ and $f_\beta(z)$ is Theorem 5.23". In the proof of this theorem (page 110), he defines $U_\beta := P_\beta / f_\beta$. Thus $U_\beta$ is a meromorphic function on the open unit disc $D(0,1)$.

Then, the author introduces some complex numbers $\omega_{j,n}\in D(0,1)$ which are simple zeros of $f_\beta$. He claims that $\omega_{j,n}$ are not poles of $U_\beta$. Note that this assertion is equivalent to say that $\omega_{j,n}$ is also a zero of $P_\beta$.
In order to show his claim, the author derives in (5.4.2) the formal identity $P_\beta(X)=U_\beta(X)f_\beta(X)$, to get 
$$
P'_\beta(X)=U'_\beta(X)f_\beta(X)+U_\beta(X)f'_\beta(X).
$$
He remarks that $f'_\beta(\omega_{j,n})\neq0$, and then {\sl specializes} $X$ to $\omega_{j,n}$. But this last manipulation is not allowed  since $U_\beta(X)$ could have a pole at $\omega_{j,n}$.\\

The same type of argument is used by the author at other places, as remarked in~\cite{Blog}.


\begin{thebibliography}{[Bo--Gi--SO]}

\bibitem{Do} E.~Dobrowolski, ``On a question of Lehmer and the number of irreducible factors of a polynomial", {\sl Acta Arith.}, {\bf 34} (1979), 391--401.

\bibitem{Le} D.H.~Lehmer, ``Factorization of certain cyclotomic functions"; Ann. of Math., {\bf 34} (1933), 461--479.

\bibitem{Blog} Math{\sl overflow}. https://mathoverflow.net/questions/286640/proof-of-the-conjecture-of-lehmer-a-dobrowolski-type-minoration

\bibitem{JLVG} J.-L. Verger-Gaugry ``A proof of the conjecture of Lehmer and of the conjecture of Schinzel-Zassenhaus".  https://arxiv.org/pdf/1709.03771.pdf and https://hal.archives-ouvertes.fr/hal-01584495v1


\end{thebibliography}
\end{document}